\documentclass[12pt]{article}
\usepackage{e-jc}

\usepackage{amsmath,amsthm,amssymb,graphicx}

\usepackage[colorlinks=true,citecolor=black,linkcolor=black,urlcolor=blue]{hyperref}

\theoremstyle{plain}
\newtheorem{theorem}{Theorem}
\newtheorem{lemma}[theorem]{Lemma}
\newtheorem{corollary}[theorem]{Corollary}
\newtheorem{prop}[theorem]{Proposition}

\theoremstyle{definition}

\theoremstyle{remark}

\begin{document}

\title{On the number of walks in a triangular domain}

\author{Paul R.~G.~Mortimer and Thomas Prellberg\\
\small School of Mathematical Sciences\\[-0.8ex]
\small Queen Mary University of London\\[-0.8ex] 
\small Mile End Road, London E1 4NS, UK\\
\small\tt \{p.r.g.mortimer,t.prellberg\}@qmul.ac.uk
}

\date{\dateline{Feb 18, 2014}{}\\
\small Mathematics Subject Classifications: 05A15, 05.10}

\maketitle

\begin{center}
  \rule{10cm}{1pt}
\end{center}

\begin{abstract}
We consider walks on a triangular domain that is a subset of the triangular lattice. We then specialise this by dividing the lattice into two directed sublattices with different weights. Our central result is an explicit formula for the generating function of walks starting at a fixed point in this domain and ending 
anywhere within the domain. Intriguingly, the specialisation of this formula to walks starting in a fixed corner
of the triangle shows that these are equinumerous to two-coloured Motzkin paths, and two-coloured three-candidate Ballot paths, in a strip of finite height.
\end{abstract}

\begin{center}
  \rule{10cm}{1pt}
\end{center}

\section{Introduction}

Recently there has been significant development of the Kernel method \cite{banderier,banderier2,prodinger}, a technique in enumerative combinatorics. This method can be used to solve linear combinatorial functional equations in so-called catalytic variables. 

While the Kernel method has been reasonably well understood when only one catalytic variable is involved, once there are two or more 
catalytic variables the situation is far from clear. There are indications that the structure of the solution, such as whether the generating function is 
algebraic or even differentiably finite, depends on the group of symmetries of the kernel of the functional equation
\cite{mishna2010,mishna}.  Only recently has there been some progress using a multivariate Kernel method in a special case \cite{bousquet2}.

Based on the Kernel method, it is possible to derive generating functions for counting problems in previously inaccessible situations. As an example, 
the exact solution of a lattice model of partially directed walks in a wedge has only been possible using an iterative 
version of the Kernel method \cite{57}, leading to a generating function that is not differentiably finite, as its singularities
accumulate at limit points.
This example also shows that as a by-product of enumerative combinatorics, deep combinatorial insight into connections between seemingly unrelated 
systems can be uncovered, leading to spin-off research in bijective combinatorics \cite{poznanovik,rubey}. 

The problems studied in this paper have strong links to directed walks; in particular the primary model of interest can also be formulated as two directed walks in a strip (cf Figs.~\ref{orthant} and \ref{strip}). We present a solution to an enumerative lattice path problem that is expressed in terms of a functional equation with three catalytic variables. We are able to solve this functional equation by virtue of the high symmetry of the kernel. 

The resulting generating function solution allows us to prove a result linking walks on a triangular domain to Motzkin paths. Finding a bijective proof of this result poses an intriguing open problem.

\section{Statement of Results}

\begin{figure}[t]
        \begin{center}
                \includegraphics[width=0.6\textwidth]{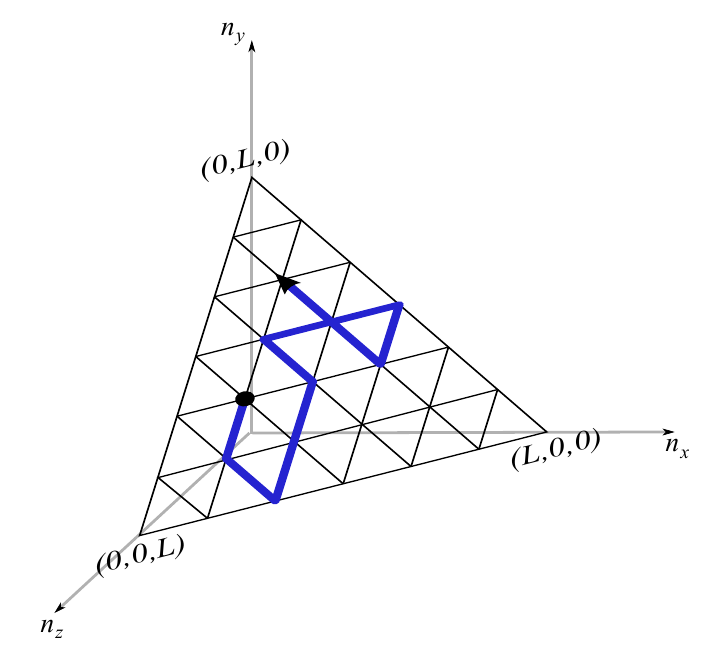}
        \end{center}
        \caption{A $10$-step walk on the bounded triangular domain intersecting the axes at $(L,0,0), (0,L,0) \text{ and } (0,0,L)$. Points on the domain are given by $(n_x,n_y,n_z)$ with $n_x,n_y,n_z\geq0$ and $n_x+n_y+n_z=L$.}
        \label{orthant}
\end{figure}

\begin{figure}[t]
        \begin{center}
                \includegraphics[width=0.6\textwidth]{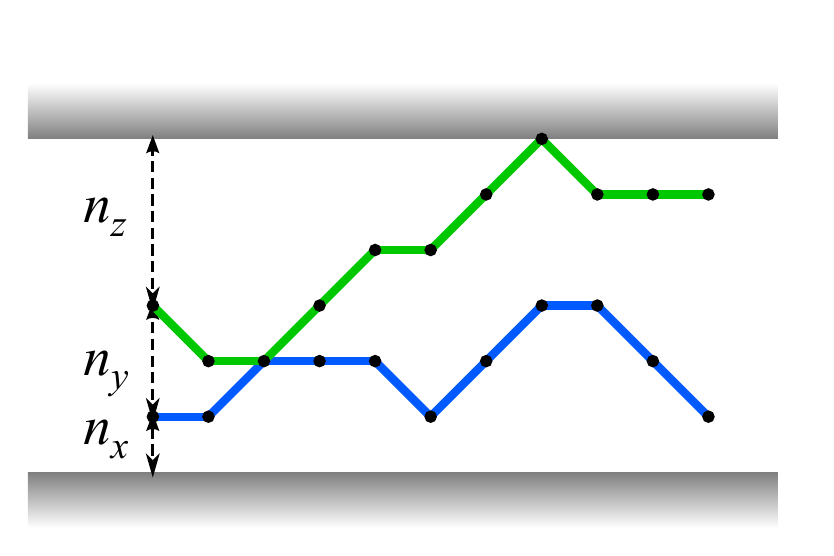}
        \end{center}
        \caption{The image of the $10$-step walk in Fig.~\ref{orthant}. The $k$-th point $(n_x,n_y,n_z)$ of the walk in Fig.~\ref{orthant} maps to two points with coordinates $(k,n_x)$ and $(k,n_x+n_y)$, respectively. This mapping generates two non-crossing directed walks that are confined to the strip $0\leq y\leq L$. In this example, the two walks touch after $2$ steps and the top walk hits the top of the strip after $7$ steps, corresponding to the points where the walk in Fig.~\ref{orthant} touches the sides of the domain. }
        \label{strip}
\end{figure}

Consider walks $(\omega_0,\omega_1,\ldots,\omega_n)$ on $\mathbb Z^{3}$ with steps $\omega_i-\omega_{i-1}$ in a step-set $\Omega_2$ such that with each step exactly one coordinate increases by one and exactly one coordinate decreases by one. More precisely, $\Omega_2= \{(1,-1,0),(-1,1,0),(1,0,-1),(-1,0,1),(0,1,-1),(0,-1,1)\}$. 

The step-set $\Omega_2$ ensures that walks lie in planes $\{(n_x,n_y,n_z)\in\mathbb Z^{3}|$ ~ $n_x+n_y+n_z=L\}$ determined by the starting point $\omega_0=(u_1,u_2,u_3)$ of the walk, where $L=u_1+u_2+u_3$. In this paper, walks on domains given by finite subsets of these planes are studied by restricting the walks to the non-negative orthant $(\mathbb N_0)^3$ (cf Fig.~\ref{orthant}). In particular, the walks lie on a bounded triangular domain; henceforth this model will be referred to as the triangle model and is of primary interest of this paper.

\begin{figure}[t]
        \begin{center}
                \includegraphics[width=0.4\textwidth]{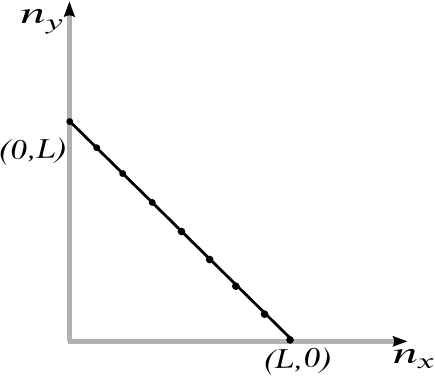}
        \end{center}
        \caption{The line model, intersecting the axes at $(L,0) \text{ and } (0,L)$. }
        \label{quadrant}
\end{figure}

As a prelude, consider walks on the $1$-dimensional analogue of the triangle model, i.e.\ walks $(\omega_0,\omega_1,\ldots,\omega_n)$ on $\mathbb Z^2$ with steps $\omega_i-\omega_{i-1}$ in a step-set $\Omega_1 = \{(1,-1),(-1,1)\}$. $\Omega_1$ ensures that walks lie on lines $\{(n_x,n_y)\in\mathbb Z^2|$ ~ $n_x+n_y=L\}$ determined by the starting point $\omega_0=(u_1,u_2)$ of the walk, where $L=u_1+u_2$. Accordingly, this will be referred to as the line model. In this paper, walks on domains given by finite subsets of these lines are studied, by restricting the walks to the non-negative quadrant $(\mathbb N_0)^2$ (cf Fig.~\ref{quadrant}).

Given a fixed starting point $\omega_0$, denote the number of $n$-step walks starting at $\omega_0$ and ending at $\omega_n=(i_1,i_2)$
by $C_n(i_1,i_2)$ and consider the generating function
\begin{equation}
G(x,y;t)=\sum_{n=0}^\infty t^n\sum_{\omega_n \in(\mathbb N_0)^2}C_n(\omega_n)x^{i_1}y^{i_2}\;,
\end{equation}
where $t$ is the generating variable conjugate to the length of the walk.
Due to the choice of the step-set $\Omega$, $G(x,y;t)$ is homogeneous of degree $L= u_1 + u_2$ in $x,y$, i.e.
\begin{equation}
G(\gamma x,\gamma y;t)=\gamma^L G(x,y;t)\;.
\end{equation}

This model has been studied before \cite[pages~7-8]{brak2005}, and has obvious connections to the generating function of Chebyshev polynomials. It is easy to solve, and gives the following result.

\begin{prop}
\label{dim1prop}
The generating function $G(x,y;t)$, which counts $n$-step walks 
starting at fixed $\omega_0=(u,v)$, is given by

\begin{equation}
G(x,y;t)=\frac 1{1-\dfrac{\frac{x}{y}+\frac{y}{x}}{p+\frac{1}{p}}}
\left(x^uy^v - \frac{x^{u+v+1}p^{v+1}(1-p^{2u+2})}{y(1-p^{2u+2v+4})} - \frac{y^{u+v+1}p^{u+1}(1-p^{2v+2})}{x(1-p^{2u+2v+4})}\right)\;,
\end{equation}
where
\begin{equation}
p=\frac{1-\sqrt{1-4t^2}}{2t}
\end{equation}
is the generating function of Dyck paths.
\end{prop}
This simplifies considerably when specifying $x=y=1$.
\begin{corollary}
\label{dim1cor}
The generating function $G(1,1;t)$, which counts $n$-step walks 
starting at fixed $\omega_0=(u,v)$ with no restrictions on the endpoint, is given by
\begin{equation}
G(1,1;t)=\frac{(1+p^2)(1-p^{u+1})(1-p^{v+1})}{(1-p)^2(1+p^{u+v+2})},
\end{equation} 
where
\begin{equation}
p=\frac{1-\sqrt{1-4t^2}}{2t}.
\end{equation}
\end{corollary}

\begin{figure}[t]
        \begin{center}
                \includegraphics[width=0.6\textwidth]{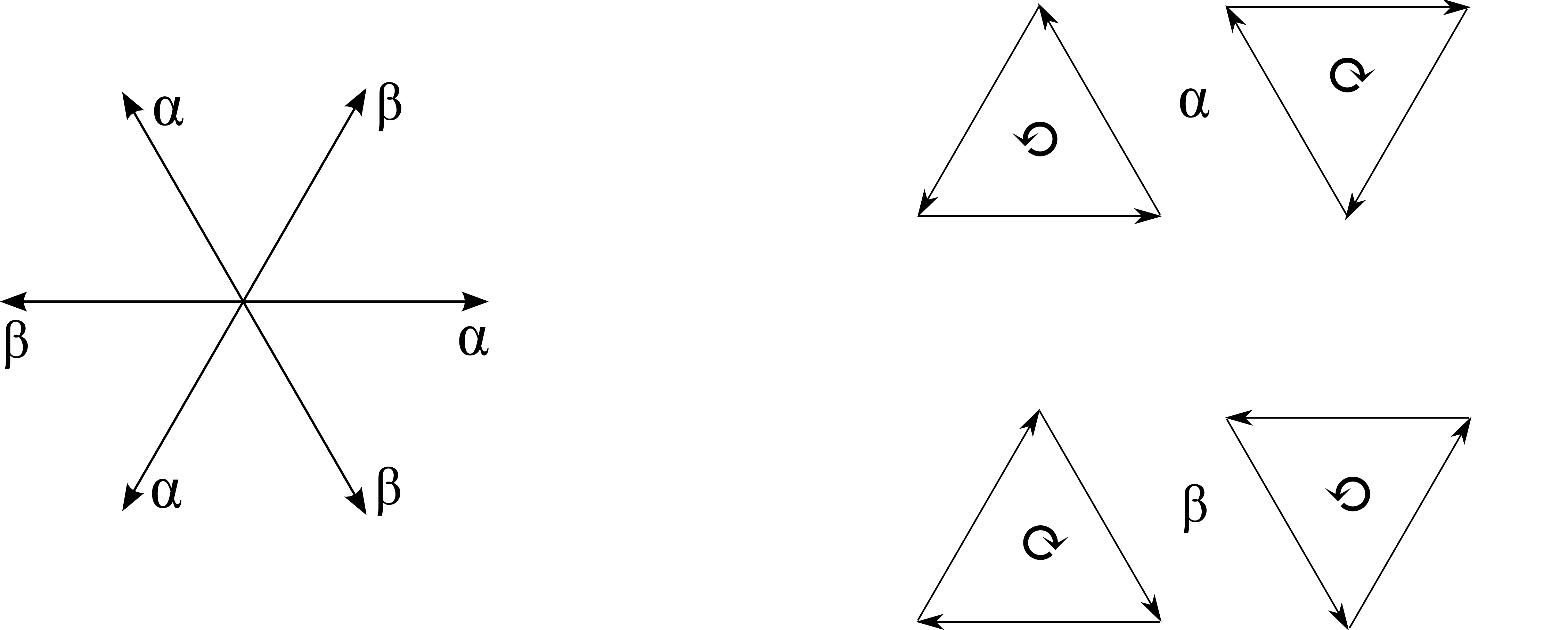}
        \end{center}
        \caption{The step-set for the triangle partitioned into two smaller step-sets $\Omega_2'$ and $\Omega_2''$, with associated weights $\alpha$ and $\beta$ respectively. Using only steps from $\Omega_2'$ restricts to allowing steps in only one orientation (clockwise or anti-clockwise) around each triangle in the domain, with neighbouring triangles (those sharing an edge) permitting opposite orientations. Using only steps from $\Omega_2''$ reverses these orientations. }
        \label{sublatticesfigure}
\end{figure}

In this paper, the main results concern the triangle model, or more precisely, a weighted generalisation of this model. Partition $\Omega_2$ into
\begin{equation}
\label{sublattices}
\begin{array}{ll}
\Omega_2'=\{(1,0,-1),(-1,1,0),(0,-1,1)\}\quad\text{and}\\
\Omega_2''=\{(1,-1,0),(-1,0,1),(0,1,-1)\}\;,
\end{array}
\end{equation}
with steps in $\Omega_2'$ and $\Omega_2''$ given the weights $\alpha$ and $\beta$, respectively (cf Fig.~\ref{sublatticesfigure}).

Given a fixed starting point $\omega_0$, denote the number of $n$-step walks starting at $\omega_0$ and ending at $\omega_n=(i_1,i_2,i_3)$
by $C_n(i_1,i_2,i_3)$ and consider the generating function
\begin{equation}
G(x,y,z;t)=\sum_{n=0}^\infty t^n\sum_{\omega_n \in(\mathbb N_0)^{d+1}}C_n(\omega_n)x^{i_1}y^{i_2}z^{i_3}\;,
\end{equation}
where $t$ is the generating variable conjugate to the length of the walk.
Due to the choice of the step-set $\Omega_2$, $G(x,y,z;t)$ is homogeneous of degree $L=u_1+u_2+u_3$ in $x,y,z$, i.e.
\begin{equation}
G(\gamma x,\gamma y,\gamma z;t)=\gamma^L G(x,y,z;t)\;.
\end{equation}
The main result of this paper is as follows.

\begin{theorem}
\label{mainthm}
The generating function $G(t)\equiv G(1,1,1;t)$, which counts $n$-step walks 
starting at fixed $\omega_0=(u,v,w)$ with no restrictions on the endpoint, is given by
\begin{equation}
\label{mainres}
G(t)=\frac{(1-p^3)(1-p^{u+1})(1-p^{v+1})(1-p^{w+1})}{(1-p)^3(1-p^{u+v+w+3})}\;,
\end{equation}
with
\begin{equation}
p=(\alpha + \beta)t M((\alpha + \beta) t)
\end{equation}
where
\begin{equation}
M(t)=\frac{1-t-\sqrt{(1+t)(1-3t)}}{2t^2}
\end{equation}
is the generating function of Motzkin paths.
\end{theorem}

For walks starting in a corner of a triangle of side-length $L$ one finds the following intriguing equinumeracy result.

\begin{corollary}
\label{equinumerous}~

\begin{itemize}
\item[(a)] Walks starting in a corner of a triangle of side-length $L=2H+1$ with arbitrary endpoint and taking $p$ steps on $\Omega_2'$ and $q$ steps on $\Omega_2''$ are in bijection with two-coloured Motzkin paths in a strip of height $H$ which have $p$ steps coloured with colour $A$ and $q$ steps coloured with colour $B$.
\item[(b)] Walks starting at a corner of a triangle of side-length $L=2H$ with arbitrary endpoint and taking $p$ steps on $\Omega_2'$ and $q$ steps on $\Omega_2''$ are in bijection with two-coloured Motzkin paths in a strip of height $H$, such that horizontal steps at height $H$ are forbidden, which have $p$ steps coloured with colour $A$ and $q$ steps coloured with colour $B$.
\end{itemize}
\end{corollary}

In particular, this immediately implies that walks starting in a corner of a triangle of side-length $L=2H+1$ (resp. $L=2H$) with arbitrary endpoint are in bijection with two-coloured Motzkin paths in a strip of height $H$ (resp. a strip of height $H$, such that horizontal steps at height $H$ are forbidden). Additionally, setting $q=0$ implies that walks starting in a corner of a triangle of side-length $L=2H+1$ (resp. $L=2H$) with arbitrary endpoint, which only take steps on $\Omega_2'$, are in bijection with Motzkin paths in a strip of height $H$ (resp. a strip of height $H$ such that horizontal steps at height $H$ are forbidden).

For walks starting in the centre of a triangle of side-length $L=3u$, there is a further result.

\begin{prop}
\label{propsymm}
The generating function $g(t)\equiv G(1,1,0;t)$, which counts walks starting at $\omega_0=(u,u,u)$ and ending at a fixed side of the triangle,
is given by
\begin{equation}
g(t)=p^u\frac{(1-p^3)(1-p^{u+1})}{(1-p)(1-p^{3u+3})}\;,
\end{equation}
with $p$ as in Theorem \ref{mainthm}.
\end{prop}

Define an $n$-step three-candidate Ballot path to be a walk $(\omega_0,\omega_1,\ldots,\omega_n)$ on $\mathbb N_0^2$ starting at the origin with steps $\omega_i-\omega_{i-1}$ taken from the step-set 
\begin{equation}
\Delta=\{(1,1),(1,-1),(1,0)\}\;,
\end{equation}
such that after $r$ steps the number of $(1,1)$ steps is greater than or equal to the number of $(1,-1)$ steps, which is greater than or equal to the number of $(1,0)$ steps, for all $0\leq r\leq n$. These can also be thought of as a coding for Yamanouchi words with three letters \cite[page~6]{azenhas2009}. Define further an $n$-step three-candidate Ballot path with excess $L$ to be an $n$-step Ballot path such that after $r$ steps the difference between the number of $(1,1)$ steps and $(1,0)$ steps is at most $L$, for all $0\leq r\leq n$.
\begin{prop}
\label{ballot}
Walks starting in a corner of a triangle of side-length $L$ with arbitrary endpoint, restricted to the sublattice $\Omega_2'$, are in bijection with three-candidate Ballot paths with excess $L$.
\end{prop}

\section{Proofs}
\subsection{Line Model}

The line model will be examined first, with the same techniques then being applied to the triangle model.

\subsubsection{Functional Equation}

An $n$-step walk is uniquely constructed by appending a step from the step-set $\Omega$ to an $(n-1)$-step walk, provided $n>0$.
This leads to the following functional equation for the generating function $G(x,y;t)$.

\begin{align}
\label{feqn1'}
G(x,y;t) = & x^{u_1}y^{u_2} + G(x,y;t)t\bigg( \frac{x}{y} + \frac{y}{x} \bigg) \nonumber \\
& - G(x, 0 ;t)t \bigg( \frac{x}{y} \bigg) -  G(0,y ;t)t \bigg( \frac{y}{x} \bigg)
\end{align}

Here, the monomial $x^{u_1}y^{u_2}$ corresponds to a zero-step walk starting (and ending) at $\omega_0=(u_1,u_2)$. The term $G(x,y;t)t\big( \frac{x}{y} + \frac{y}{x} \big)$ corresponds to appending any of the steps in $\Omega_1$ irrespective of whether the resulting walk steps violates the boundary condition and leaves the domain. This overcounting is adjusted by the remaining terms. For example, $G(x, 0;t)$ corresponds to walks which end at $(i_1, 0 ;t)$, and therefore $G(x, 0 ;t) \frac{x}{y}$ corresponds to precisely those walks stepping across the boundary.

As this is a functional equation for the generating function $G(x,y;t)$ in the variables $x,y$ only, the $t$-dependence is dropped
by writing. The functional equation (\ref{feqn1'}) is rewritten as
\begin{equation}
\label{feqn1}
G(x,y)\left[1-t\,\left(\frac {x}{y}+\frac {y}{x}\right)\right]=x^{u_1}y^{u_2}
-G(0,y)\,t\,\left(\frac {y}{x}\right)
-G(x,0)\,t\,\left(\frac { x}{y}\right)\;.
\end{equation}

\subsubsection{The Kernel}

Furthermore, the ``Kernel'' $K(x,y;t)\equiv K(x,y)$  of the functional equation,
\begin{equation}
\label{kernel1}
K(x,y)=1-t\left(\frac {x}{y}+\frac {y}{x}\right),
\end{equation}
is introduced. Symmetry properties of this Kernel are central to the arguments. Note that the Kernel is homogeneous of degree zero, {\em i.e.}\ it is invariant under rescaling of all the variables. This trivial symmetry will be implicitly assumed in the considerations below.

Now introduce $\mathrm G(S)$, the group of transformations which leaves the Kernel of the functional equation invariant for the step set $S$. This is in line with the notation introduced by Fayolle {\em et al.}~\cite{fayolle1999}. For the line model, the step-set is $$S_1=\left\{\frac{x}{y},\frac{y}{x}\right\},$$
where here and henceforth, steps are identified with their associated combinatorial weights.

\begin{lemma}
The Kernel $K(x,y)$ is invariant under action of the group of transformations

$$ \mathrm G(S_1)=\left \langle(y,x)\right \rangle \cong C_2\;.$$

\end{lemma}
\noindent where $(y,x)$ is shorthand notation for the map that sends $(x,y)$ to $(y,x)$.
\noindent In particular, one arrives at the following result.

\begin{lemma}
\label{invariantkernel}
The Kernel $K(x,y)$ is invariant under the following $1$-parameter substitutions.
$$K(p,1)=K(1,p)=1-t\left(p+\frac1p\right)$$
\end{lemma}
The dependence between $p$ and $t$ is henceforth fixed such that 
\begin{equation}
\label{tsubs1}
1-t(p+1/p)=0,
\end{equation}
which when solved for $p$ gives $p=t D(t)$, where $D(t)$ is the generating function for Dyck paths
\begin{equation}
\label{Dyck}
D(t)=\frac{1-\sqrt{1-4t^2}}{2t^2
}.
\end{equation}
In particular, $p$ is a well-defined power series with zero constant.

Using this dependency and substituting the two choices from Lemma $\ref{invariantkernel}$ into the functional equation (\ref{feqn1}) implies
\begin{subequations}
\label{eqns1}
\begin{align}
\label{eq_p1}tpG(p,0)+\frac{t}{p} G(0,1)&=p^u\\
\label{eq_1p}\frac tpG(1,0)+tpG(0,p)&=p^v  .
\end{align}
\end{subequations}
Using homogeneity of the generating function, replace
\begin{equation}
\label{homogeneous}
G(p,0)=p^{u+v}G(1,0)\;,\quad G(0,p)=p^{u+v}G(0,1)\;,
\end{equation}
and solve the two equations \ref{eq_p1} and \ref{eq_1p} in the two variables $G(1,0)$ and $G(0,1)$ to find that

\begin{equation}
\label{variables}
G(1,0)=\frac{p^{v+1}(p^{2u+2}-1)}{t(p^{2u+2v+4}-1)}\;,\quad
G(0,1)=\frac{p^{u+1}(p^{2v+2}-1)}{t(p^{2u+2v+4}-1)} \ .
\end{equation}
Applying the homogeneity argument of (\ref{homogeneous}) to (\ref{variables}), one determines $G(x,0)$ and $G(0,y)$, and substituting these into (\ref{feqn1}) and eliminating $t$ via (\ref{tsubs1}) gives the result stated in Proposition \ref{dim1prop}

\begin{equation}
\label{dim1result}
G(x,y;t)=\frac 1{1-\dfrac{\frac{x}{y}+\frac{y}{x}}{p+\frac{1}{p}}}
\left(x^uy^v - \frac{x^{u+v+1}p^{v+1}(1-p^{2u+2})}{y(1-p^{2u+2v+4})} - \frac{y^{u+v+1}p^{u+1}(1-p^{2v+2})}{x(1-p^{2u+2v+4})}\right)\ .
\end{equation}
Substituting $x=y=1$ into (\ref{dim1result}) then gives Corollary \ref{dim1cor}

\begin{equation}
G(1,1;t)=\frac{(1+p^2)(1-p^{u+1})(1-p^{v+1})}{(1-p)^2(1+p^{u+v+2})} \ .
\end{equation}

\subsection{Triangle Model}

The same method is now applied to the triangle model, including the weights $\alpha$ and $\beta$ corresponding to the directed sublattices (equation (\ref{sublattices})). 

\subsubsection{Functional Equation}

Again, an $n$-step walk is uniquely constructed by appending a step from the step-set $\Omega$ to an $(n-1)$-step walk, provided $n>0$. This leads to the following functional equation for the generating function $G(x,y,z;t)$.

\begin{align}
\label{feqn2'}
G(x,y,z;t) = & x^{u_1}y^{u_2}z^{u_3} + G(x,y,z;t)t\left(\frac {\beta x}{y}+\frac {\alpha y}{x}+\frac {\alpha x}{z}+\frac {\beta z}{x}+\frac {\beta y}{z}+\frac {\alpha z}{y}\right) \nonumber \\
& -G(0,y,z)\,t\,\left(\frac {\alpha y}{x}+\frac {\beta z}{x}\right)
-G(x,0,z)\,t\,\left(\frac {\beta x}{y}+\frac {\alpha z}{y}\right)\nonumber \\
& -G(x,y,0)\,t\,\left(\frac {\alpha x}{z}+\frac {\beta y}{z}\right)
\end{align}

Similarly to equation (\ref{feqn1'}), the monomial $x^{u_1}y^{u_2}z^{u_3}$ corresponds to a zero-step walk starting (and ending) at $\omega_0=(u_1,u_2,u_3)$, the second term corresponds to appending any of the steps in $\Omega$ (irrespective of whether the resulting walk steps violates the boundary condition and leaves the domain), and any overcounting is adjusted by the remaining three terms, each of which accounts for stepping over one of the three boundary edges of the triangle.

Again, as this is a functional equation for the generating function $G(x,y,z;t)$ in the variables $x,y,z$ only, the $t$-dependence is dropped
by writing $G(x,y,z;t)\equiv G(x,y,z)$. The functional equation (\ref{feqn2'}) is rewritten as
\begin{multline}
\label{feqn2}
G(x,y,z)\left[1-t\,\left(\frac {\beta x}{y}+\frac {\alpha y}{x}+\frac {\alpha x}{z}+\frac {\beta z}{x}+\frac {\beta y}{z}+\frac {\alpha z}{y}\right)\right]=x^u_1y^u_2z^u_3\\
-G(0,y,z)\,t\,\left(\frac {\alpha y}{x}+\frac {\beta z}{x}\right)
-G(x,0,z)\,t\,\left(\frac {\beta x}{y}+\frac {\alpha z}{y}\right)
-G(x,y,0)\,t\,\left(\frac {\alpha x}{z}+\frac {\beta y}{z}\right)\;.
\end{multline}

\subsubsection{The Kernel}

Again, the Kernel $K(x,y,z;t)\equiv K(x,y,z)$  of the functional equation,
\begin{equation}
\label{kernel2}
K(x,y,z)=1-t\left(\frac {\beta x}{y}+\frac {\alpha y}{x}+\frac {\alpha x}{z}+\frac {\beta z}{x}+\frac {\beta y}{z}+\frac {\alpha z}{y}\right)\;,
\end{equation}
is needed, and symmetry properties of this Kernel are central to the arguments. As before, note that the Kernel is homogeneous of degree zero, {\em i.e.}\ it is invariant under rescaling of all the variables, and this trivial symmetry will be implicitly assumed in the considerations below.

Looking again at $\mathrm G(S)$, the group of transformations which leaves the Kernel of the functional equation invariant for the step set $S$, the step-set for the triangle model is
$$S_2=\left\{\frac {\beta x}{y}+\frac {\alpha y}{x}+\frac {\alpha x}{z}+\frac {\beta z}{x}+\frac {\beta y}{z}+\frac {\alpha z}{y}\right\}\;,$$
and $ \mathrm G(S_2)$ is generated by a rotation and an inversion.
\begin{lemma}
The Kernel $K(x,y,z)$ is invariant under action of the group of transformations

$$ \mathrm G(S_2)=\left\langle (y,z,x),\left(\frac{1}{x},\frac{1}{y},\frac{1}{z}\right)\right\rangle \cong C_3 \times C_2\;.$$
\end{lemma}

\noindent Moreover, there is a one-variable
sub-set which has useful consequences.
\begin{lemma}
\label{thereshouldbe}
The Kernel $K(x,y,z)$ is invariant under the following $1$-parameter substitutions.
\begin{multline}
K(1,1,p)=K(1,p,1)=K(p,1,1)=K(1,p,p)=K(p,1,p)=K(p,p,1)\\
=1-t(\alpha + \beta)\left(p+1+\frac1p\right)\;.
\end{multline}
\end{lemma}

\noindent Fixing the dependence between $p$ and $t$ such that 
\begin{equation}
\label{tsubs2}
1-t(\alpha + \beta)(p+1+1/p)=0\;,
\end{equation}
gives $p=(\alpha + \beta)t M((\alpha + \beta) t)$, where
\begin{equation}
M(t)=\frac{1-t-\sqrt{(1+t)(1-3t)}}{2t^2}
\end{equation}
is the generating function of Motzkin paths. In particular, $p$ is a well-defined power series in $t$ with zero constant.

Using this dependency and substituting the six choices from Lemma \ref{thereshouldbe} into the functional equation (\ref{feqn2}) then implies
\begin{subequations}
\label{eqns2}
\begin{align}
\label{eq_p11}\frac{(\alpha + \beta)t}pG(0,1,1)+t(\alpha +\beta p)G(p,0,1)+t(\alpha p +\beta )G(p,1,0)&=p^u\\
\label{eq_1p1}\frac{(\alpha + \beta)t}pG(1,0,1)+t(\alpha p +\beta)G(0,p,1)+t(\alpha +\beta p)G(1,p,0)&=p^v\\
\label{eq_11p}\frac{(\alpha + \beta)t}pG(1,1,0)+t(\alpha +\beta p)G(0,1,p)+t(\alpha p + \beta)G(1,0,p)&=p^w\\
\label{eq_1pp}(\alpha + \beta)tpG(0,p,p)+t\left(\alpha +\frac{\beta}{p}\right)G(1,0,p)+t\left(\frac{\alpha}{p} + \beta \right)G(1,p,0)&=p^vp^w\\
\label{eq_p1p}(\alpha + \beta)tpG(p,0,p)+t\left(\frac{\alpha}{p} + \beta \right)G(0,1,p)+t\left(\alpha +\frac{\beta}{p}\right)G(p,1,0)&=p^up^w\\
\label{eq_pp1}(\alpha + \beta)tpG(p,p,0)+t\left(\alpha +\frac{\beta}{p}\right)G(0,p,1)+t\left(\frac{\alpha}{p} + \beta \right)G(p,0,1)&=p^up^v
\end{align}
\end{subequations}
Using homogeneity of the generating function, replace
\begin{equation}
G(p,p,0)=p^LG(1,1,0)\;,\quad G(p,0,p)=p^LG(1,0,1)\;,\quad G(0,p,p)=p^LG(0,1,1)\;,
\end{equation}
and from the linear combination $[(\ref{eq_p11})+(\ref{eq_1p1})+(\ref{eq_11p})]-p[(\ref{eq_1pp})+(\ref{eq_p1p})+(\ref{eq_pp1})]$ it is easily found that
\begin{equation}
\label{linearcomb}
(\alpha + \beta)t[G(0,1,1)+G(1,0,1)+G(1,1,0)]=
\frac{p^{u+1}+p^{v+1}+p^{w+1}-p^{2+L}(p^{-u}+p^{-v}+p^{-w})}{1-p^{3+L}}\;.
\end{equation}
Substituting $(x,y,z)=(1,1,1)$ into (\ref{feqn2}) shows that $G(1,1,1)$ can be computed explicitly, as
\begin{equation}
\label{G111}
(1-3(\alpha + \beta)t)G(1,1,1)=1-(\alpha + \beta)t[G(0,1,1)+G(1,0,1)+G(1,1,0)]\;.
\end{equation}
Substituting (\ref{linearcomb}) into (\ref{G111}) and eliminating $t$ via (\ref{tsubs2}) gives the desired
final result
\begin{equation}
G(1,1,1)=\frac{(1-p^3)(1-p^{u+1})(1-p^{v+1})(1-p^{w+1})}{(1-p)^3(1-p^{3+L})}\;.
\end{equation}
Finally, note that substituting $p(t)=(\alpha + \beta)tM((\alpha + \beta)t)$ followed by $(\alpha + \beta)t=s$ into (\ref{tsubs2}) implies
\begin{equation}
M(s)=1+sM(s)+s^2M(s)^2\;,
\end{equation}
whence $M(s)$ is the Motzkin path generating function. This completes the proof of Theorem \ref{mainthm}.

Letting $(u,v,w)=(L,0,0)$ in (\ref{mainres}) implies that the generating function for walks starting in a
corner is given by
\begin{equation}
\label{thisformula}
G(1,1,1)=\frac{(1-p^3)(1-p^{1+L})}{(1-p)(1-p^{3+L})}\;.
\end{equation}

\subsubsection{Continued Fractions}

Equation (\ref{thisformula}) is intimately related to the convergents of the continued fraction expansion of the Motzkin path generating function. One can show by mathematical induction that in this case $G(1,1,1)$ can be written as a continued fraction. More precisely, for $L=2H$ even, there is a continued fraction of length $H$,
\begin{align}
\label{Leven}
\underbrace{\cfrac1{1-(\alpha + \beta)t-\cfrac{(\alpha + \beta)^2t^2}{1-(\alpha + \beta)t-\cfrac{(\alpha + \beta)^2t^2}{\ddots-\cfrac{(\alpha + \beta)^2t^2}{1-(\alpha + \beta)t-(\alpha + \beta)^2t^2}}}}}_{\text{length $H$}}&=\frac{(1-p^3)(1-p^{1+2H})}{(1-p)(1-p^{3+2H})}\;,
\end{align}
and for $L=2H+1$ odd, there is a continued fraction of length $H+1$,
\begin{align}
\label{Lodd}
\underbrace{\cfrac1{1-(\alpha + \beta)t-\cfrac{(\alpha + \beta)^2t^2}{1-(\alpha + \beta)t-\cfrac{(\alpha + \beta)^2t^2}{\ddots-\cfrac{(\alpha + \beta)^2t^2}{1-(\alpha + \beta)t}}}}}_{\text{length $H+1$}}&=\frac{(1-p^3)(1-p^{2+2H})}{(1-p)(1-p^{4+2H})}\;.
\end{align}
It is easy to show that equations (\ref{Leven}) and (\ref{Lodd}) hold for the base case $H=0$,
\begin{subequations}
\begin{align}
1&=\frac{(1-p^3)(1-p^{1+0})}{(1-p^{1+0})(1-p^{3+0})}\\
\cfrac1{1-(\alpha + \beta) t}&=\frac{(1-p^3)(1-p^{2+0})}{(1-p)(1-p^{4+0})}\;,
\end{align}
\end{subequations}
and the inductive step follows from showing that
\begin{equation}
\frac{(1-p^3)(1-p^{1+(L+2)})}{(1-p)(1-p^{3+(L+2)})}=\cfrac1{1-(\alpha+\beta)t-(\alpha+\beta)^2t^2\cfrac{(1-p^3)(1-p^{1+L})}{(1-p)(1-p^{3+L})}}\;.
\end{equation}

From the combinatorial theory of continued fractions the combinatorial interpretation in terms of Motzkin paths follows easily, as given in a paper by Flajolet \cite[pages~6-11]{flajolet1980}. This immediately implies Corollary \ref{equinumerous}. Substituting $\alpha = 1, \ \beta=1$  into (\ref{Leven}) and (\ref{Lodd}) gives coefficients $2$ and $4$, of $t$ and $t^2$ respectively, following from the fact that the relevant Motzkin paths are two-coloured, and substituting in $\alpha = 1, \ \beta=0$ gives the interpretation in terms of normal Motzkin paths.

\subsubsection{Further Results}

Attempting to solve the triangle model in full generality proved beyond the reach of the techniques used in this paper, as the system of equations (\ref{eqns2}) is underdetermined, linking nine quantities with six equations. The only case in which one can extract further information from it is one of high symmetry, namely when the starting point is chosen to be in the centre of the triangle, i.e. $\omega_0=(u,u,u)$, in which the triangle has size $L=3u$. The equations (\ref{eqns2}) then reduce to two equations in two unknowns,
\begin{subequations}
\label{symmeqns}
\begin{align}
\label{symmeq_p11}\frac{(\alpha + \beta)t}pG(1,1,0)+(\alpha +\beta)t(1+p)G(p,1,0)&=p^u\\
\label{symmeq_1pp}(\alpha + \beta)tp^{1+3u}G(1,1,0)+(\alpha + \beta)t\left(1+\frac1p\right)G(p,1,0)&=p^{2u}\;.
\end{align}
\end{subequations}
This can readily be solved, and
\begin{equation}
(\alpha + \beta)tG(1,1,0)=
\frac{p^{1+u}(1-p^{1+u})}{1-p^{3+3u}}\;.
\end{equation}
Eliminating $t$ by using (\ref{tsubs2}) proves Proposition \ref{propsymm}.

It now remains to prove our final result. Without loss of generality, starting in the corner marked by coordinates $(L,0,0)$, the steps $(-1,0,1)$, $(0,1,-1)$, and $(1,-1,0)$ can be mapped to $(1,1)$, $(1,-1)$, and $(1,0)$, respectively. This maps steps in $\omega'$ to steps in three-candidate Ballot paths and the restrictions imposed by the boundaries of the triangle clearly transfer to the restrictions on a Ballot path with excess $L$. This proves Proposition \ref{ballot}.

\section{Conclusion and Open Problems}

\subsection{The Problem in General Dimension}

We now frame the triangle model as the $2$-dimensional case of a larger class of models. Consider walks $(\omega_0,\omega_1,\ldots,\omega_n)$ on $\mathbb Z^{d+1}$ with steps $\omega_i-\omega_{i-1}$ in a step-set $\Omega_d$ such that with each step exactly one coordinate increases by one and exactly one coordinate decreases by one. More precisely, $\Omega_d$ is the set of steps with coordinates $(e_1,e_2,\hdots ,e_{d+1})$ such that for all ordered pairs $(i,j)$ with $1 \leq i,j \leq d+1$ and $i \neq j$, $e_i=1, \ e_j=-1 \text{ and } e_k=0 \text{ for all } 1 \leq k \leq d+1  \text{ and }k\notin \{i,j\}$.

The step-set $\Omega_d$ ensures that walks lie in a $d$-dimensional hyperplane $\{(n_{x_1},\hdots,n_{x_{d+1}})\in\mathbb Z^{d+1}|$ ~ $n_{x_1}+\hdots+n_{x_{d+1}}=L\}$ determined by the starting point $\omega_0=(u_1,\hdots,u_{d+1})$ of the walk, where $L=\sum_{j=1}^{d+1}u_j$. In this paper, walks on domains given by finite subsets of these hyperplanes are studied by restricting the walks to the non-negative orthant $(\mathbb N_0)^{d+1}$. 
Fixing the dimension $d$, this class of walks is referred to as the $d$-dimensional case. The $1$-dimensional case is the line model, and the $2$-dimensional case is the triangle model. In the $3$-dimensional case, the domains would be tetrahedra of side-length $L$.

Given a fixed starting point $\omega_0$, denote the number of $n$-step walks starting at $\omega_0$ and ending at $\omega_n=(i_1,\hdots,i_{d+1})$
by $C_n(i_1,\hdots,i_{d+1})$ and consider the generating function
\begin{equation}
G(x_1,\hdots ,x_{d+1};t)=\sum_{n=0}^\infty t^n\sum_{\omega_n \in(\mathbb N_0)^{d+1}}C_n(\omega_n)\prod_{j=1}^{j=d+1}x_j^{i_j}\;,
\end{equation}
where $t$ is the generating variable conjugate to the length of the walk.
Due to the choice of the step-set $\Omega_d$, $G(x_1,\hdots ,x_{d+1};t)$ is homogeneous of degree $L=\sum_{j=1}^{d+1}u_j$ in $x_1,\hdots,x_{d+1}$, i.e.
\begin{equation}
G(\gamma x_1,\hdots,\gamma x_{d+1};t)=\gamma^L G(x_1,\hdots,x_{d+1};t)\;.
\end{equation}

In this paper we have set up a general problem which we have completely solved in dimension $1$. We have also solved it in dimension $2$ in the case where endpoints are not weighted, along with a high-symmetry case. Unfortunately, our argument does not provide enough information to solve the general case for dimension $2$ case or that for higher dimensions and we must leave these open.

\subsection{Generating Function Properties}

Generating functions for walks in finite domains are rational. This is a direct result of the fact that the adjacency matrix for such systems has finite dimension. In particular, in the triangle model, a triangle of side-length $L$ contains $\binom{L+2}2$ vertices, or states, and therefore we would expect the degree of the numerator and denominator of the generating function to grow quadratically in $L$. However, for the cases where we are able to prove results, there is some cancellation such that the growth is linear in $L$. 

The process of finding our results began with some initial series generation which allowed us to predict the form of the generating functions. Using this for the general case in dimension $2$, of walks with arbitrary fixed start and end points, we have numerical evidence that in general the degrees of the numerator and denominator grow quadratically in $L$, and it may be this extra complexity that has prevented us from solving this case with our method.

\subsection{Bijections}

We have also proven an intriguing equinumeracy result in Corollary \ref{equinumerous}. Of particular interest are the two special cases noted beneath it; the equinumeracy between walks on the undirected triangular domain and two-coloured Motzkin paths, and that between walks on one of the directed sublattices and normal Motzkin paths.  Taking only steps on the directed sublattice $\Omega'$ halves the out-degree of every vertex in the domain, and so it is clear that the result for two-coloured Motzkin paths implies the result for normal Motzkin paths.

If a bijective proof of Corollary \ref{equinumerous} were to be found, this might elucidate the connections between this model and continued fractions, and thus open avenues towards solving other models. Eu \cite{eu2010} gives a bijective proof of the directed case for triangular domains of infinite side-length via standard Young tableaux (which are a coding of Yamanouchi words), Yeats \cite{yeats} gives a bijective proof of the undirected case for domains of infinite side-length using intermediate markings, and the authors of this paper have bijective proofs of the undirected case for side-lengths $L=1,2 \text{ and } 3$. We note that Proposition \ref{ballot} provides a possible alternative route to a bijective proof, via three-candidate Ballot paths. However, we have not been able to find a proof for general finite side-length, and therefore leave this as an open problem.

\end{document}